\documentclass{article}
\usepackage{amsfonts}\usepackage{amssymb}
\usepackage{graphicx}\usepackage{amsmath}
\usepackage{abstract}
\usepackage{color}
\usepackage{graphics}   \usepackage{amsthm}
\newtheorem{theorem}{Theorem}[section]
\newtheorem{lemma}[theorem]{Lemma}

\numberwithin{equation}{section}
\providecommand{\keywords}[1]{\textbf{Key words.} #1}

\title{NEW ERROR ESTIMATES OF LINEAR TRIANGLE FINITE ELEMENTS FOR THE STEKLOV EIGENVALUE PROBLEM}

\author{Hai Bi, Yidu Yang, Yuanyuan Yu
 \\\\
{\small School of Mathematical Sciences, Guizhou Normal University,}\\{\small  Guiyang,  $550001$,  China}
\\{\small
bihaimath@gznu.edu.cn, ydyang@gznu.edu.cn, yuyuanyuan567@126.com}
}

\begin{document}
\date{}
\maketitle

\begin{abstract}
In this paper we make a further discussion on the finite elements approximation for the Steklov
eigenvalue problem on concave polygonal domain. We make full use of the regularity estimate and
the characteristic of edge average interpolation operator of nonconforming Crouzeix-Raviart element, which is different from the existing proof argument, and prove a new and optimal error estimate in $\|\cdot\|_{0,\partial\Omega}$ for the eigenfunction of linear conforming finite element and the nonconforming Crouzeix-Raviart element, which is an improvement of the current results. Finally, we present some numerical experiments to support the theoretical analysis.
\end{abstract}


\section{Introduction}

Steklov eigenvalue problems have important physical
background and many applications. For instance, they appear
in the analysis of stability of mechanical oscillators immersed
in a viscous fluid (see \cite{conca} and the references therein), in the study
of surface waves (see \cite{bergman}), in the study of the vibration modes
of a structure in contact with an incompressible fluid (see
\cite{bermudez}), in the analysis of the antiplane shearing on a system of collinear
faults under slip-dependent friction law (see \cite{bucur}), etc. Thus
the numerical methods for solving these problems have attracted more and more scholars'
attention. Till now, systematical and profound studies on the
conforming finite elements approximation for Steklov eigenvalue
problems have been made on polygonal domain such as
\cite{andreev,armentano1,armentano2,bermudez,bramble,garau,lim,xie,yang7}).
Recently, the nonconforming finite elements for Steklov problems have
also been considered, e.g., see \cite{alonso,bi,liq2,liu,yang1}. The aim of this paper is to discuss the error estimates of linear triangle finite elements, including the linear conforming finite element and the nonconforming Crouzeix-Raviart element, approximation for Steklov eigenvalue problems with variable coefficients on concave polygonal domain.

We consider the following Steklov eigenvalue problem
\begin{eqnarray}\label{e1.1}
-div(\alpha\nabla u) +\beta u=0~~in~ \Omega,~~~~
\alpha\frac{\partial u}{\partial n}=\lambda u~~ on~ \partial\Omega,
\end{eqnarray}
where $\Omega \subset R^{2}$ is a polygonal domain with $\omega$
being the largest inner angle of $\Omega$, and
$\frac{\partial u}{\partial n}$ is the outward normal derivative.

Having in mind that $H^{s}(\Omega)$ denotes the Sobolev
space with real order $s$ on $\Omega$, $\|\cdot\|_{s}$ is the norm
on $H^{s}(\Omega)$ and $H^{0}(\Omega)=L_{2}(\Omega)$, and
$H^{s}(\partial\Omega)$ denotes the Sobolev space with real order
$s$ on $\partial\Omega$ with the norm $\|\cdot\|_{s,\partial\Omega}$.

Suppose that the coefficients $\alpha=\alpha(x)$ and
$\beta=\beta(x)$ are bounded by above and below by positive
constants. We assume that $\alpha\in C^{1}(\bar{\Omega})$.

The weak form of (\ref{e1.1}) is given by: Find $\lambda\in
R$, $u \in H^{1}(\Omega),\|u\|_{0,\partial\Omega}=1$, such that
\begin{eqnarray}\label{e1.2}
a(u,v)=\lambda b(u,v),~~~\forall v\in H^{1}(\Omega),
\end{eqnarray}
where
\begin{eqnarray*}
a(u,v)=\int\limits_{\Omega}(\alpha\nabla u\cdot \nabla v+\beta uv)
dx,~b(u,v)=\int \limits_{\partial\Omega}uvds.
\end{eqnarray*}
It is easy to know that $a(\cdot,\cdot)$ is a symmetric, continuous
and $H^{1}(\Omega)$-elliptic bilinear form
on $H^{1}(\Omega)\times H^{1}(\Omega)$.

In the existing literatures, the error estimate of linear
triangle elements eigenfunction, including conforming element and nonconforming
Crouzeix-Raviart element (hereafter termed C-R element for simplicity), in
$\|\cdot\|_{0, \partial\Omega}$ is all $O(h^{r+\frac{r}{2}})$ where
$r$ is the regularity exponent of the eigenfunction (see Lemma 2.1).
It is obvious that this estimate is not optimal since it doesn't achieve the order of
interpolation error. In this paper, we improve this estimate when eigenfunctions are
singular (i.e., $r<1$) and prove that in this case the error
estimate of linear triangle elements eigenfunction can achieve
$O(h^{r+\frac{1}{2}})$. Comparing the proof arguments of existing estimates (e.g., see \cite{armentano1,bramble,liq2,yang1}), we make full use of the regularity estimate and the characteristic of edge average interpolation operator of C-R element, especially in the analysis for conforming finite elements, and obtain the improved error estimates (\ref{e2.25}) and (\ref{e3.6}) which are optimal.

Throughout this paper, $C$ denotes a positive constant
independent of $h$, which may not be the same constant in different
places.

\section{The nonconforming Crouzeix-Raviart element approximation for the Steklov eigenvalue problem}

Consider the source problem (\ref{e2.1}) associated with (\ref{e1.1}): Find $w \in H^{1}(\Omega)$, such that
\begin{eqnarray}\label{e2.1}
a(w,v)=b(f,v),~~~\forall v \in H^{1}(\Omega).
\end{eqnarray}
As for the source problem  (\ref{e2.1}), there hold the following regularity results.
\begin{lemma}
If $f\in L_{2}(\partial\Omega)$,
then $w\in H^{1+\frac{r}{2}}(\Omega)$ and
\begin{eqnarray}\label{e2.2}
\|w\|_{1+\frac{r}{2}}\leq C_{\Omega}\|f\|_{0,\partial\Omega};
\end{eqnarray}
if $f\in H^{\frac{1}{2}}(\partial\Omega)$, then $w\in
H^{1+r}(\Omega)$ and
\begin{eqnarray}\label{e2.3} \|w\|_{1+r}\leq
C_{\Omega}\|f\|_{\frac{1}{2},\partial\Omega};
\end{eqnarray}
if $f\in H^{\varepsilon}(\partial\Omega)$,
$\varepsilon\in(0,r-1/2)$, then $w\in H^{\frac{3}{2}+\varepsilon}(\Omega)$
and
\begin{eqnarray}\label{e2.4}
\|w\|_{\frac{3}{2}+\varepsilon}\leq C_{\Omega}\|f\|_{\varepsilon,\partial\Omega}.
\end{eqnarray}
Here $r=1$ when $\omega<\pi$, and $r<\frac{\pi}{\omega}$ which can
be arbitrarily close to $\frac{\pi}{\omega}$ when $\omega>\pi$, and $C_{\Omega}$ is a priori constant.
\end{lemma}

{\it Proof.} See \cite{dauge}. ~~~$\Box$

Note that $a(\cdot,\cdot)$ is coercive, using the source
problem (\ref{e2.1}) associated with (\ref{e1.2}) we can define the
operator $A: L_{2}(\partial\Omega)\to H^{\frac{3}{2}}(\Omega)\subset
H^{1}(\Omega)$, satisfying
\begin{eqnarray*}
a(Af,v)&=&b(f,v),~~~ \forall v \in H^{1}(\Omega).
\end{eqnarray*}
Define the operator $T: L_{2}(\partial\Omega)\to H^{1}(\partial\Omega)$ satisfying
\begin{eqnarray*}
Tf=(Af)',
\end{eqnarray*}
where $'$ denotes the restriction to $\partial\Omega$.

Bramble and Osborn \cite{bramble} proved that (\ref{e1.2})
has the operator form:
\begin{eqnarray}\label{e2.5}
Tw=\mu w.
\end{eqnarray}
Namely, if $(\mu,w)$ is an eigenpair of (\ref{e2.5}), then
$(\lambda,Aw)$ is an eigenpair of (\ref{e1.2}),
$\lambda=\frac{1}{\mu}$. Conversely, if $(\lambda, u)$ is an
eigenpair of (\ref{e1.2}), then $(\mu,u')$ is an eigenpair of
(\ref{e2.5}),
$\mu=\frac{1}{\lambda}$.

Let $\lambda$ be the $j$-th eigenvalue of $T$. We
arrange eigenvalues by the increasing order with each eigenvalue
counted according to its algebraic multiplicity. And let
$M(\lambda)$ denote the space spanned by eigenfunctions of
(\ref{e1.2}) corresponding to the eigenvalue $\lambda$.

Let $\pi_{h}=\{K\}$ be a regular triangulation of $\Omega$
in the sense of the minimal internal angle condition (see
\cite{ciarlet}, pp. 131). We denote $h=\max_{K\in \pi_{h}}h_{K}$
where $h_{K}$ is the diameter of element $K$. Let $S^{h}_{nc}$
be the C-R element space (see \cite{crouzeix}) defined on $\pi_{h}$:\\
\indent $S^{h}_{nc}=\{v\in L_{2}(\Omega):v\mid_{K} \in
span\{1,x_{1},x_{2}\}$, $v$ is
continuous at the midpoints of the edges of elements$\}$.

The C-R element approximation of
(\ref{e1.2}) is: Find $\lambda_{h}\in R$, $u_{h} \in
S^{h}_{nc},\|u_{h}\|_{0,\partial\Omega}=1$, such that
\begin{eqnarray}\label{e2.6}
a_{h}(u_{h},v)=\lambda_{h} b(u_{h},v),~~~\forall v\in S^{h}_{nc},
\end{eqnarray}
where
\begin{eqnarray*}
a_{h}(u_{h},v)=\sum\limits_{K\in\pi_{h}}\int\limits_{K}(\alpha\nabla
u_{h}\cdot \nabla v+\beta u_{h}v) dx.
\end{eqnarray*}
\indent Define $\|v\|_{h}=(\sum\limits_{K\in
\pi_{h}}\|v\|_{1,K}^{2})^{\frac{1}{2}}$,
$\|v\|_{1,K}^{2}=\int_{K}(|\frac{\partial v }{\partial
x_{1}}|^{2}+|\frac{\partial v }{\partial x_{2}}|^{2}+|v|^{2})dx$.
Evidently, $\|\cdot\|_{h}$ is the norm on $S^{h}_{nc}$ and it is simple
to show that $a_{h}(\cdot,\cdot)$ is uniformly $S^{h}_{nc}$-elliptic.

The C-R element approximation of (\ref{e2.1})
is: Find $w_{h} \in S^{h}_{nc}$, such that
\begin{eqnarray}\label{e2.7}
a_{h}(w_{h},v)=b(f,v),~~~\forall v \in S^{h}_{nc}.
\end{eqnarray}

Denote the consistency term of the C-R element by
\begin{eqnarray}\label{e2.8}
E_{h}(w,v)=a_{h}(w,v)-b(f,v).
\end{eqnarray}
And based on the standard method (see, for example
\cite{alonso,bi,liq2}), the following consistency error estimate can be
proved.

\begin{lemma}
Suppose that $w\in H^{1+t}(\Omega)$ with
$t\in [\frac{r}{2}, 1]$ is the weak solution of (\ref{e2.1}), then
\begin{eqnarray}\label{e2.9}
&&E_{h}(w,v)\leq C h^{t}\|w\|_{1+t}\|v\|_{h},~~~\forall v\in
S^{h}_{nc}+H^{1}(\Omega),\\\label{e2.10} &&\|w-w_{h}\|_{h}\leq
Ch^{t}\|w\|_{1+t},\\\label{e2.11} &&\|w-w_{h}\|_{0,\partial\Omega}\leq
Ch^{t+\frac{r}{2}}\|w\|_{1+t}.
\end{eqnarray}
\end{lemma}

{\it Proof.} See, e.g., Lemma 2.2 in \cite{bi}, or Theorem 2.1 in \cite{liq2}.~~~ $\Box$

Define the interpolation operator $I_{h}:H^{1}(\Omega)\to
S^{h}_{nc}$:
\begin{eqnarray}\label{e2.12}
\int_{l}I_{h}uds=\int_{l}uds~~~\forall l,~\forall u \in
H^{1}(\Omega),
\end{eqnarray}
where $l$ is an edge of arbitrary element in $\pi_{h}$.

According to the interpolation theory (see \cite{ciarlet}), we have
\begin{eqnarray}\label{e2.13}
\|u-I_{h}u\|_{0}&\leq& Ch^{1+r}\|u\|_{1+r},\\\label{e2.14}
\|u-I_{h}u\|_{h}&\leq& Ch^{r}\|u\|_{1+r}.
\end{eqnarray}

\begin{theorem}
Let $w\in H^{1+t}(\Omega)$ with $t\in
(\frac{1}{2}, 1]$ be the weak solution of (\ref{e2.1}), then
\begin{eqnarray}\label{e2.15}
\|w-w_{h}\|_{-\varepsilon, \partial\Omega}\leq
Ch^{t+\frac{1}{2}+\varepsilon}\|w\|_{1+t}
\end{eqnarray}
where $t=\frac{1}{2}+\varepsilon$ if $f\in
H^{\varepsilon}(\partial\Omega)$, $0<\varepsilon<r-\frac{1}{2}$, and
$t=r$ if $f\in H^{\frac{1}{2}}(\partial\Omega)$.
\end{theorem}

{\it Proof.} For each $g\in
H^{\varepsilon}(\partial\Omega)$, let $\varphi$ be the unique
solution of the following variational problem:
\begin{eqnarray*}
a(v,\varphi)=b(g,v),~~\forall v\in H^{1}(\Omega).
\end{eqnarray*}
From (\ref{e2.4}) we know that $\varphi \in
H^{\frac{3}{2}+\varepsilon}(\Omega)$. Let
$\varphi_{h}=I_{h}\varphi\in S^{h}$ be the interpolation of
$\varphi$, then
\begin{eqnarray*}
&&b(f,\varphi-\varphi_{h})=b(f,\varphi)-b(f,\varphi_{h})\\
&&~~~=a_{h}(w,\varphi)-a_{h}(w_{h},\varphi_{h})\\
&&~~~=a_{h}(w,\varphi)-a_{h}(w_{h},\varphi)+a_{h}(w_{h},\varphi)-a_{h}(w_{h},\varphi_{h})\\
&&~~~=a_{h}(w-w_{h},\varphi)+a_{h}(w_{h}-w,\varphi-\varphi_{h})+a_{h}(w,\varphi-\varphi_{h}).
\end{eqnarray*}
By the definition of consistency term we have
\begin{eqnarray*}
&&a_{h}(w-w_{h},\varphi)+a_{h}(w,\varphi-\varphi_{h})\\
&&~~=a_{h}(w-w_{h},\varphi)-b(w-w_{h},g)+b(w-w_{h},g)\\
&&~~~+a_{h}(w,\varphi-\varphi_{h})-b(f,\varphi-\varphi_{h})+b(f,\varphi-\varphi_{h})\\
&&~~=E_{h}(\varphi,w-w_{h})+b(w-w_{h},g)+E_{h}(w,\varphi-\varphi_{h})+b(f,\varphi-\varphi_{h}).
\end{eqnarray*}
Combining the above two relationships, we get
\begin{eqnarray*}
b(w-w_{h},g)=-E_{h}(\varphi,w-w_{h})-E_{h}(w,\varphi-\varphi_{h})
-a_{h}(w_{h}-w,\varphi-\varphi_{h}),
\end{eqnarray*}
then
\begin{eqnarray}\label{e2.16}
&&|b(g,w-w_{h})|=|b(w-w_{h},g)|\nonumber\\
&&~~~\leq |E_{h}(\varphi,w-w_{h})|+|E_{h}(w,\varphi-\varphi_{h})|
+|a_{h}(w_{h}-w,\varphi-\varphi_{h})|.
\end{eqnarray}
From (\ref{e2.9}), (\ref{e2.3}), (\ref{e2.4}), (\ref{e2.10}) and the
error estimate of interpolation, we can deduce that
\begin{eqnarray*}
&&|E_{h}(\varphi,w-w_{h})|\leq
Ch^{\frac{1}{2}+\varepsilon}\|\varphi\|_{\frac{3}{2}+\varepsilon}\|w-w_{h}\|_{h}
\leq C h^{t+\frac{1}{2}+\varepsilon}\|w\|_{1+t}\|g\|_{\varepsilon,\partial\Omega},\\
&&|E_{h}(w,\varphi-\varphi_{h})|\leq
Ch^{t}\|w\|_{1+t}\|\varphi-\varphi_{h}\|_{h}
\leq C h^{t+\frac{1}{2}+\varepsilon}\|w\|_{1+t}\|g\|_{\varepsilon,\partial\Omega},\\
&&|a_{h}(w_{h}-w,\varphi-\varphi_{h})|\leq
C\|w-w_{h}\|_{h}\|\varphi-\varphi_{h}\|_{h}\leq
Ch^{t+\frac{1}{2}+\varepsilon}\|w\|_{1+t}\|g\|_{\varepsilon,\partial\Omega}.
\end{eqnarray*}
And substituting the above three estimates into (\ref{e2.16}), we
obtain
\begin{eqnarray*}
|b(g,w-w_{h})|\leq
Ch^{t+\frac{1}{2}+\varepsilon}\|w\|_{1+t}\|g\|_{\varepsilon,\partial\Omega}.
\end{eqnarray*}
By the definition of negative norm, we have
\begin{eqnarray}\label{e2.17}
&&\|w-w_{h}\|_{-\varepsilon,\partial\Omega}
=\sup\limits_{g\in H^{\varepsilon}(\partial\Omega)}\frac{|b(g,w-w_{h})|}{\|g\|_{\varepsilon,\partial\Omega}}\nonumber\\
&&~~~\leq \sup\limits_{g\in
H^{\varepsilon}(\partial\Omega)}\frac{Ch^{t+\frac{1}{2}+\varepsilon}\|w\|_{1+t}\|g\|_{\varepsilon,\partial\Omega}}{\|g\|_{\varepsilon,\partial\Omega}}
\leq C h^{t+\frac{1}{2}+\varepsilon}\|w\|_{1+t},
\end{eqnarray}
namely, (\ref{e2.15}) is true.~~~$\Box$

Let us denote by $\delta S^{h}_{nc}$ the functions defined on
$\partial\Omega$, which are restriction of functions in $S^{h}_{nc}$ to
$\partial\Omega$. From \cite{alonso}, pp.189 we know that
\begin{eqnarray}\label{e2.18}
\delta S^{h}\subset H^{\varepsilon}(\partial\Omega),~~~ for~any~\varepsilon \in [0,\frac{1}{2}).
\end{eqnarray}

Since $a_{h}(\cdot,\cdot)$ is uniformly elliptic with
respect to $h$, the approximate source problem (\ref{e2.7})
associated with (\ref{e2.6}) is uniquely solvable.
 Thus, we can define the
operator $A_{h}:L_{2}(\partial\Omega)\to S^{h}_{nc}$, satisfying
\begin{eqnarray*}
a_{h}(A_{h}f,v)&=&b(f,v),~~~ \forall v \in S^{h}_{nc}.
\end{eqnarray*}
Define $T_{h}: L_{2}(\partial\Omega)\to \delta S^{h}_{nc}\subset
H^{\varepsilon}(\partial\Omega)\subset L_{2}(\partial\Omega)$, satisfying
\begin{eqnarray*}
T_{h}f=(A_{h}f)'.
\end{eqnarray*}
\indent \cite{yang1} proved that (\ref{e2.6}) has the operator form:
\begin{eqnarray}\label{e2.19}
T_{h}w_{h}=\mu_{h} w_{h}.
\end{eqnarray}
Namely, if $(\mu_{h},w_{h})$ is an eigenpair of (\ref{e2.19}), then
$(\lambda_{h}, A_{h}w_{h})$ is an eigenpair of (\ref{e2.6}),
$\lambda_{h}=\frac{1}{\mu_{h}}$. Conversely, if $(\lambda_{h},
u_{h})$ is an eigenpair of (\ref{e2.6}), then $(\mu_{h},u_{h}')$ is
an eigenpair of (\ref{e2.19}),
$\mu_{h}=\frac{1}{\lambda_{h}}$.

We prove the following interpolation estimates.
\begin{lemma}
Let $u\in H^{1+r}(\Omega)$, then the
following estimates hold:
\begin{eqnarray}\label{e2.20}
&&||u-I_{h}u||_{0, \partial\Omega}\leq C
h^{r+\frac{1}{2}}\|u\|_{1+r},\\\label{e2.21}
&&||u-I_{h}u||_{-\varepsilon, \partial\Omega}\leq C
h^{r+\frac{1}{2}+\varepsilon}\|u\|_{1+r}.
\end{eqnarray}
\end{lemma}

{\it Proof.}  Let $l\subset\partial\Omega$ be the edge of
the element $K$, then by the trace inequality ( see Lemma 7.1.1 in
\cite{wang}) we have
\begin{eqnarray*}
||u-I_{h}u||_{0, l}\leq C
(h_{K}^{-\frac{1}{2}}\|u-I_{h}u\|_{0,K}+h_{K}^{\frac{1}{2}}\|u-I_{h}u\|_{1,K})\leq
Ch_{K}^{r+\frac{1}{2}}\|u\|_{1+r, K},
\end{eqnarray*}
thus (\ref{e2.20}) is valid.\\
For any $g\in H^{\varepsilon}(\partial\Omega)~(\varepsilon\in
(0,r-\frac{1}{2}))$, let $I_{0}g$ be the piecewise constant
interpolation of $g$ on $\partial\Omega$. From the definition of
$I_{h}$ and the interpolation estimates we have
\begin{eqnarray*}
|b(g,u-I_{h}u)|&=&|b(g-I_{0}g,u-I_{h}u)|\\
&\leq& C \|g-I_{0}g\|_{0, \partial\Omega}\|u-I_{h}u\|_{0, \partial\Omega}\\
&\leq& Ch^{r+\frac{1}{2}+\varepsilon}\|u\|_{1+r}\|g\|_{\varepsilon,\partial\Omega},
\end{eqnarray*}
and using the definition of
negative norm we know that (\ref{e2.21}) holds.~~~$\Box$

\begin{lemma}
Suppose that $M(\lambda)\subset
H^{1+r}(\Omega)$ and $\lambda$ be the $j$-th eigenvalue of
(\ref{e1.2}). Let $\lambda_{h}$  be the $j$-th eigenvalue of
(\ref{e2.6}) and $u_{h}$ be an eigenfunction corresponding to
$\lambda_{h}$ with $\|u_{h}\|_{0,\partial\Omega}=1$. Then there exists $u\in
M(\lambda)$ with $\|u\|_{0,\partial\Omega}=1$, such that
\begin{eqnarray}\label{e2.22}
&&\mid\lambda_{h}-\lambda\mid \leq C h^{2r}\|
u\|_{1+r}^{2},\\\label{e2.23}
&&\|u-u_{h}\|_{h}\leq C  h^{r}\|
u\|_{1+r},\\\label{e2.24}
&&\|u-u_{h}\|_{0,\partial\Omega}\leq C h^{\frac{3}{2}r}\|
u\|_{1+r}.
\end{eqnarray}
\end{lemma}

{\it Proof.} See \cite{alonso,liq2,yang1}.~~~$\Box$

Lemma 2.5 is an existing conclusion. Next we will improve the estimate (\ref{e2.24}).

\begin{theorem}
Under the conditions of Lemma 2.5,
further assume that $\pi_{h}$ is a quasi-uniform mesh (see pp.135 in
\cite{ciarlet}), then
\begin{eqnarray}\label{e2.25}
\|u-u_{h}\|_{0,\partial\Omega}&\leq& C h^{r+\frac{1}{2}}\|
u\|_{1+r}.
\end{eqnarray}
\end{theorem}

{\it Proof.} Since $Au$ and $A_{h}u$ are solutions of (\ref{e2.1}) and
(\ref{e2.7}) with $f=u$, respectively, then from (\ref{e2.15}) we
know that
\begin{eqnarray}\label{e2.26}
&~&\|Tu-T_{h}u\|_{-\varepsilon,\partial\Omega}\leq
Ch^{r+\frac{1}{2}+\varepsilon}\|
u\|_{1+r}.
\end{eqnarray}
Using (\ref{e2.21}) we obtain
\begin{eqnarray}\label{e2.27}
\|Tu-I_{h}Tu\|_{-\varepsilon,
\partial\Omega}\leq C h^{r+\frac{1}{2}+\varepsilon}\|u\|_{1+r}.
\end{eqnarray}
From (\ref{e2.26}) and (\ref{e2.27}), we have
\begin{eqnarray}\label{e2.28}
\|T_{h}u-I_{h}Tu\|_{-\varepsilon,
\partial\Omega}\leq C h^{r+\frac{1}{2}+\varepsilon}\|u\|_{1+r}.
\end{eqnarray}
By the definition of negative norm and the inverse estimates, we
have
\begin{eqnarray*}
&&\|T_{h}u-I_{h}Tu\|_{0,
\partial\Omega}^{2}
\leq \|T_{h}u-I_{h}Tu\|_{-\varepsilon,
\partial\Omega}\|T_{h}u-I_{h}Tu\|_{\varepsilon,
\partial\Omega}\\
&&\leq \|T_{h}u-I_{h}Tu\|_{-\varepsilon,
\partial\Omega}h^{-\varepsilon}\|T_{h}u-I_{h}Tu\|_{0,
\partial\Omega},
\end{eqnarray*}
thus
\begin{eqnarray}\label{e2.29}
&&\|T_{h}u-I_{h}Tu\|_{0,
\partial\Omega}
\leq C h^{r+\frac{1}{2}}\|u\|_{1+r}.
\end{eqnarray}
By using (\ref{e2.29}) and (\ref{e2.20}), we get
\begin{eqnarray}\label{e2.30}
\|Tu-T_{h}u\|_{0,
\partial\Omega}
&\leq & \|T_{h}u-I_{h}Tu\|_{0,
\partial\Omega}+\|Tu-I_{h}Tu\|_{0,
\partial\Omega}\nonumber\\
&\leq& Ch^{r+\frac{1}{2}}\|u\|_{1+r}.
\end{eqnarray}
It has been proved in \cite{liq2,bi} that
$\|T-T_{h}\|_{0,\partial\Omega}\to 0(h\to 0)$, thus, from Theorem
7.4 in \cite{babuska} we get
\begin{eqnarray}\label{e2.31}
\|u_{h}-u\|_{0,\partial\Omega}\leq
C\|(T-T_{h})u\|_{0,\partial\Omega}.
\end{eqnarray}
Substituting (\ref{e2.30}) into (\ref{e2.31}), we
obtain (\ref{e2.25}).~~~ $\Box$

\noindent{\bf Remark 2.1.} If $r\in (\frac{1}{2},1)$, i.e.,
$\Omega$ is concave, it is clear that the estimate (\ref{e2.25}) is
better than (\ref{e2.24}).

\section{The conforming element approximation for the Steklov eigenvalue problem}

Let $S^{h}_{c}\subset C(\overline{\Omega})$ be a space of piecewise linear polynomials defined on $\pi_{h}$. The conforming
element approximation of (\ref{e1.2}) is: Find $\lambda_{h}\in R$,
$u_{h} \in S^{h}_{c}$ with $\|u_{h}\|_{0,\partial\Omega}=1$, such that
\begin{eqnarray}\label{e3.1}
a(u_{h},v)=\lambda_{h} b(u_{h},v),~~~\forall v\in S^{h}_{c}.
\end{eqnarray}

As for the conforming finite element approximation (\ref{e3.1}), the following results are valid (see \cite{armentano1,bramble}).\\
\begin{lemma}
Let $(\lambda_{h}, u_{h})$ be the $j$-th
eigenpair of
 (\ref{e3.1}), $\lambda$ be the $j$-th
eigenvalue of
 (\ref{e1.2}), and $M(\lambda)\subset H^{1+r}(\Omega)$.
Then there exists $u\in M(\lambda)$ such that
\begin{eqnarray}\label{e3.2}
 &&|\lambda-\lambda_{h}|
\leq C h^{2r}\|u\|_{1+r},\\\label{e3.3}
 &&\|u-u_{h}\|_{1} \leq C
h^{r}\|u\|_{1+r},
\\\label{e3.4} &&\|u-u_{h}\|_{0,\partial\Omega}
\leq C h^{\frac{3r}{2}}\|u\|_{1+r},
\end{eqnarray}
where the principle to determine $r$ see Lemma 2.1.
\end{lemma}

Now, let $P_{h}: H^{1}(\Omega)\to
 S^{h}_{c}$ be the Ritz projection defined by
\begin{eqnarray*}
a(w-P_{h}w, v)=0,~~~\forall v\in S^{h}_{c}.
\end{eqnarray*}
We can define the operator $A_{h}:L_{2}(\partial\Omega)\to S^{h}_{c}$,
satisfying
\begin{eqnarray*}
a(A_{h}f,v)&=&b(f,v),~~~ \forall v \in S^{h}_{c}.
\end{eqnarray*}
It is easy to know that $A_{h}=P_{h}A$.

Let $\delta S^{h}_{c}$ be the space of functions defined on
$\partial\Omega$, which are restriction of functions in $S^{h}_{c}$ to
$\partial\Omega$. Define $T_{h}: L_{2}(\partial\Omega)\to \delta S^{h}_{c}\subset
H^{1}(\partial\Omega)$, satisfying
\begin{eqnarray*}
T_{h}f=(A_{h}f)'.
\end{eqnarray*}
It has been proved in \cite{babuska,bramble} that
$\|T-T_{h}\|_{0,\partial\Omega}\to 0~(h\to 0)$, and (\ref{e3.1}) has
the operator form:
\begin{eqnarray}\label{e3.5}
T_{h}w_{h}=\frac{1}{\lambda_{h}} w_{h}.
\end{eqnarray}

Next we will give a new error estimate for the conforming finite element.
\begin{theorem}
Under the conditions of Lemma 3.1,
further assume that $\pi_{h}$ is quasi-uniform mesh, then
\begin{eqnarray}\label{e3.6}
\|u-u_{h}\|_{0,\partial\Omega}\leq C h^{r+\frac{1}{2}}\| u\|_{1+r}.
\end{eqnarray}
\end{theorem}

{\it Proof.} For each $g\in
H^{\varepsilon}(\partial\Omega)$, let $\varphi$ be the unique
solution of the following variational problem:
\begin{eqnarray*}
a(v,\varphi)=b(g,v),~~\forall v\in H^{1}(\Omega).
\end{eqnarray*}
From (\ref{e2.4}) we know that $Ag=\varphi \in
H^{\frac{3}{2}+\varepsilon}(\Omega)$, and
\begin{eqnarray*}
&&b(u-P_{h}u,g)=a(u-P_{h}u,Ag)\\
&&~~~=a(u-P_{h}u,Ag-P_{h}Ag)\leq
Ch^{r+\frac{1}{2}+\varepsilon}\|u\|_{1+r}\|g\|_{\varepsilon, \partial\Omega},
\end{eqnarray*}
thus, by the definition of negative norm, we have
\begin{eqnarray}\label{e3.7}
&&\|u-P_{h}u\|_{-\varepsilon, \partial\Omega}\leq
Ch^{r+\frac{1}{2}+\varepsilon}\|u\|_{1+r}.
\end{eqnarray}
Let $I_{h}u$ be the interpolation of $u$ defined by (\ref{e2.12}).
By using the inverse estimates, (\ref{e3.7}), (\ref{e2.21}) and
(\ref{e2.20}), we get
\begin{eqnarray}\label{e3.8}
&&\|u-P_{h}u\|_{0,
\partial\Omega}
\leq \|P_{h}u-I_{h}u\|_{0,
\partial\Omega}+\|u-I_{h}u\|_{0,
\partial\Omega}\nonumber\\
&&~~~\leq  Ch^{-\varepsilon}
\|P_{h}u-I_{h}u\|_{-\varepsilon,
\partial\Omega}+\|u-I_{h}u\|_{0,
\partial\Omega}\nonumber\\
&&~~~\leq  Ch^{-\varepsilon}
(\|P_{h}u-u\|_{-\varepsilon,
\partial\Omega}+\|u-I_{h}u\|_{-\varepsilon,
\partial\Omega})+\|u-I_{h}u\|_{0,
\partial\Omega}\nonumber\\
&&~~~\leq C h^{r+\frac{1}{2}}\|u\|_{1+r}.
\end{eqnarray}
By using the spectral approximation theory, we get
\begin{eqnarray}\label{e3.9}
\|u-u_{h}\|_{0,\partial\Omega}&\leq& C\|Tu-T_{h}u\|_{0,\partial\Omega}=
C\|Au-A_{h}u\|_{0,\partial\Omega}\nonumber\\
&=&C\|Au-P_{h}Au\|_{0,\partial\Omega}
= C\frac{1}{\lambda}\|u-P_{h}u\|_{0,\partial\Omega}.
\end{eqnarray}
Substituting (\ref{e3.8}) into (\ref{e3.9}), we obtain
(\ref{e3.6}).~~~$\Box$
\\

\noindent{\bf Remark 3.1.} Comparing (\ref{e3.4}) and (\ref{e3.6}), we can see that when eigenfunctions
are singular, i.e., $r<1$, the error estimate in
$\|\cdot\|_{0,\partial\Omega}$ is improved.\\
\indent When we prove the improved estimates (\ref{e2.25}) and (\ref{e3.6}), we make full use of the
regularity estimate (\ref{e2.4}) to analyze the negative norm estimate, then use the negative norm
estimate and the interpolation of C-R element, especially in the analysis for conforming elements,
to obtain the optimal estimates in $L^{2}(\partial\Omega)$; while the existing work is to analyze directly the error in $L^{2}(\partial\Omega)$ by using (\ref{e2.2}) which leads to the lost of error order.

\noindent{\bf Remark 3.2.} We prove the estimates (\ref{e2.25}) and (\ref{e3.6}) under the condition that
$\pi_{h}$ is quasi-uniform. In fact, this condition is not a restriction. Since
when $\tilde{\pi}_{h}$ is a regular partition derived from $\pi_{h}$ by local refinement,
the approximate eigenfunction $\tilde{u}_{h}$ computed on $\tilde{\pi}_{h}$ generally satisfies
$\|\tilde{u}_{h}-u\|_{0,\partial\Omega}\leq C \|u_{h}-u\|_{0,\partial\Omega}$,
then, for such regular meshes (\ref{e2.25}) and (\ref{e3.6}) are still valid.

\section{Numerical Experiments}

Consider the problem (\ref{e1.1}), where
$\alpha(x)=\beta(x)=1$, $\Omega \subset R^{2} $, $\Omega=[0,1]\times
[0,\frac{1}{2}]\bigcup [0,\frac{1}{2}]\times [\frac{1}{2},1]$ is a
L-shaped domain with the largest inner angle
$\omega=\frac{3}{2}\pi$, or $\Omega=[0,1]\times [0,1]\setminus
\{x=(x_{1},x_{2}):\frac{1}{2}\leq x_{1}\leq 1, x_{2}=\frac{1}{2}\}$
is the unit square with a slit which the largest inner angle $\omega=2\pi$.

We adopt a uniform isosceles right triangulation $\pi_{h}$.
We use the formula $ratio(\lambda_{h})=lg(\frac{\lambda_{h}-\lambda}{\lambda_{h/2}-\lambda})/lg2$
and $ratio(u_{h})=lg(\frac{\|u_{h}-u\|_{0,\partial\Omega}}{\|u_{h/2}-u\|_{0,\partial\Omega}})/lg2$
to compute the convergence order of approximations of linear conforming element to validate our analysis.

By calculation we find that the eigenfunction associated with
$\lambda_{2}$ is singular. So in our numerical experiments we compute the approximation of the second
eigenvalue $\lambda_{2,h}$ and the corresponding eigenfunction $u_{2,h}$.
Since the exact eigenpairs of the problem (\ref{e1.1})
are unknown, we use the adaptive method to compute a high-precision approximation $\lambda_{2}\in [0.89364476, 0.89364690]$ for the L-shaped domain and $\lambda_{2}\in [0.734554376, 0.73455822]$ for the unit square with a slit, and use them as the exact values, and the corresponding eigenfunction $u$ is taken as the approximation computed on the uniform mesh with the mesh diameter $h=\frac{\sqrt{2}}{1024}$. The numerical results on the L-shaped domain and the slit domain are listed in Table 1 and Table 2, respectively.

\begin{center}{\bf Table 1: The results by using linear conforming element on the
L-shaped domain}\end{center}
\begin{center}\footnotesize
\begin{tabular}{ccccc}\hline
$h$&$\lambda_{2,h}$&$ratio(\lambda_{2,h})$&$\|u_{2,h}-u\|_{0,\partial\Omega}$&$ratio(u_{2,h})$\\
\hline
$\frac{\sqrt{2}}{8}$ &0.92115806 &1.40979290  &0.02800065 & 1.12103345\\
$\frac{\sqrt{2}}{16}$ &0.90400049&1.39401631 &0.01287370 & 1.12243165 \\
$\frac{\sqrt{2}}{32}$ &0.89758582&1.37720866 &0.00591313 & 1.13940548\\
$\frac{\sqrt{2}}{64}$&0.89516258 & 1.36395810 &0.00268425 & 1.16362010\\
$\frac{\sqrt{2}}{128}$&0.89423511& 1.35431173 &0.00119822 & 1.20908566\\
$\frac{\sqrt{2}}{256}$&0.89387631&  &0.00051828 & \\\hline

\end{tabular}
\end{center}

\begin{center}{\bf Table 2: The results by using linear conforming element on the
unit square with a slit}\end{center}
\begin{center}\footnotesize
\begin{tabular}{ccccc}\hline
$h$&$\lambda_{2,h}$&$ratio(\lambda_{2,h})$&$\|u_{2,h}-u\|_{0,\partial\Omega}$&$ratio(u_{2,h})$\\
\hline
$\frac{\sqrt{2}}{8}$ &0.79372467&  1.05162089 &0.04741663 & 0.85372994\\
$\frac{\sqrt{2}}{16}$ &0.76310065& 1.03053027 &0.02623810 & 0.88264565\\
$\frac{\sqrt{2}}{32}$ &0.74852962& 1.01703653 &0.01423081 & 0.91530246\\
$\frac{\sqrt{2}}{64}$&0.74146094&  1.00934037 &0.00754564 & 0.95510109\\
$\frac{\sqrt{2}}{128}$&0.73798634 &1.00505527 &0.00389208  & 1.01983389\\
$\frac{\sqrt{2}}{256}$&0.73626532 & &0.00191947  & \\\hline
\end{tabular}
\end{center}

For the L-shaped domain $r=\frac{2}{3}$, $2r=\frac{4}{3}$. From Table 1 we can see that the convergence order of $\lambda_{2,h}$ is approximately equal to
$2r=\frac{4}{3}\approx 1.333333$. It also can be seen from Table 1 that the convergence order of
$u_{2,h}$ is very close to $r+\frac{1}{2}=\frac{7}{6}\approx 1.166667$, which is coincide with
the theoretical result (\ref{e3.6}); while the convergence order of $u_{2,h}$ according to the previous conclusion (\ref{e3.4}) should be $\frac{3r}{2}=1$.

For the unit square with a slit $r=\frac{1}{2}$. From Table 2 we can see that the convergence order of $\lambda_{2,h}$ is approximately equal to
$2r=1$. We can also see from Table 2 that the convergence order of
$u_{2,h}$ is very close to $r+\frac{1}{2}=1$, which is coincide with
the theoretical result (\ref{e3.6}); while the previous conclusion (\ref{e3.4}) states that the convergence order of $u_{2,h}$ is $\frac{3r}{2}=0.75$.\\

\noindent {\bf Acknowledgments.}  This work was supported by the
National Natural Science Foundation of China (Grant Nos. 11201093, 10761003).

\end{document}